\documentclass{amsart} %article
\usepackage{amsfonts,amssymb,amsmath}%,euscript}
    \newtheorem{thm}{Theorem}                     [section]
    \newtheorem{thm*}{Theorem}
    \newtheorem{prop}[thm]{Proposition}
    \newtheorem{lemma}[thm]{Lemma}
    \newtheorem{cor}[thm]{Corollary}

    \newtheorem{lemma*}{Lemma}    %for AMS \newtheorem*
    \newtheorem{rems*}{Remark}   %for AMS \newtheorem*
\newcommand{\ndef}{\newcommand*}
\def\rndef{\renewcommand}

\ndef{\myaddress}[1]{\begin{center} \it\small #1 \end{center}}

\ndef{\clA}{{\mathcal A}} \ndef{\rmA}{{\mathrm A}} \ndef{\mbA}{{\mathbb A}} \ndef{\bfA}{{\mathbf A}} \ndef{\euA}{{\EuScript A}} \ndef{\frA}{{\mathfrak A}}
\ndef{\clB}{{\mathcal B}} \ndef{\rmB}{{\mathrm B}} \ndef{\mbB}{{\mathbb B}} \ndef{\bfB}{{\mathbf B}} \ndef{\euB}{{\EuScript B}} \ndef{\frB}{{\mathfrak B}}
\ndef{\clC}{{\mathcal C}} \ndef{\rmC}{{\mathrm C}} \ndef{\mbC}{{\mathbb C}} \ndef{\bfC}{{\mathbf C}} \ndef{\euC}{{\EuScript C}} \ndef{\frC}{{\mathfrak C}}
\ndef{\clD}{{\mathcal D}} \ndef{\rmD}{{\mathrm D}} \ndef{\mbD}{{\mathbb D}} \ndef{\bfD}{{\mathbf D}} \ndef{\euD}{{\EuScript D}} \ndef{\frD}{{\mathfrak D}}
\ndef{\clE}{{\mathcal E}} \ndef{\rmE}{{\mathrm E}} \ndef{\mbE}{{\mathbb E}} \ndef{\bfE}{{\mathbf E}} \ndef{\euE}{{\EuScript E}} \ndef{\frE}{{\mathfrak E}}
\ndef{\clF}{{\mathcal F}} \ndef{\rmF}{{\mathrm F}} \ndef{\mbF}{{\mathbb F}} \ndef{\bfF}{{\mathbf F}} \ndef{\euF}{{\EuScript F}} \ndef{\frF}{{\mathfrak F}}
\ndef{\clG}{{\mathcal G}} \ndef{\rmG}{{\mathrm G}} \ndef{\mbG}{{\mathbb G}} \ndef{\bfG}{{\mathbf G}} \ndef{\euG}{{\EuScript G}} \ndef{\frG}{{\mathfrak G}}
\ndef{\clH}{{\mathcal H}} \ndef{\rmH}{{\mathrm H}} \ndef{\mbH}{{\mathbb H}} \ndef{\bfH}{{\mathbf H}} \ndef{\euH}{{\EuScript H}} \ndef{\frH}{{\mathfrak H}}
\ndef{\clI}{{\mathcal I}} \ndef{\rmI}{{\mathrm I}} \ndef{\mbI}{{\mathbb I}} \ndef{\bfI}{{\mathbf I}} \ndef{\euI}{{\EuScript I}} \ndef{\frI}{{\mathfrak I}}
\ndef{\clJ}{{\mathcal J}} \ndef{\rmJ}{{\mathrm J}} \ndef{\mbJ}{{\mathbb J}} \ndef{\bfJ}{{\mathbf J}} \ndef{\euJ}{{\EuScript J}} \ndef{\frJ}{{\mathfrak J}}
\ndef{\clK}{{\mathcal K}} \ndef{\rmK}{{\mathrm K}} \ndef{\mbK}{{\mathbb K}} \ndef{\bfK}{{\mathbf K}} \ndef{\euK}{{\EuScript K}} \ndef{\frK}{{\mathfrak K}}
\ndef{\clL}{{\mathcal L}} \ndef{\rmL}{{\mathrm L}} \ndef{\mbL}{{\mathbb L}} \ndef{\bfL}{{\mathbf L}} \ndef{\euL}{{\EuScript L}} \ndef{\frL}{{\mathfrak L}}
\ndef{\clM}{{\mathcal M}} \ndef{\rmM}{{\mathrm M}} \ndef{\mbM}{{\mathbb M}} \ndef{\bfM}{{\mathbf M}} \ndef{\euM}{{\EuScript M}} \ndef{\frM}{{\mathfrak M}}
\ndef{\clN}{{\mathcal N}} \ndef{\rmN}{{\mathrm N}} \ndef{\mbN}{{\mathbb N}} \ndef{\bfN}{{\mathbf N}} \ndef{\euN}{{\EuScript N}} \ndef{\frN}{{\mathfrak N}}
\ndef{\clO}{{\mathcal O}} \ndef{\rmO}{{\mathrm O}} \ndef{\mbO}{{\mathbb O}} \ndef{\bfO}{{\mathbf O}} \ndef{\euO}{{\EuScript O}} \ndef{\frO}{{\mathfrak O}}
\ndef{\clP}{{\mathcal P}} \ndef{\rmP}{{\mathrm P}} \ndef{\mbP}{{\mathbb P}} \ndef{\bfP}{{\mathbf P}} \ndef{\euP}{{\EuScript P}} \ndef{\frP}{{\mathfrak P}}
\ndef{\clQ}{{\mathcal Q}} \ndef{\rmQ}{{\mathrm Q}} \ndef{\mbQ}{{\mathbb Q}} \ndef{\bfQ}{{\mathbf Q}} \ndef{\euQ}{{\EuScript Q}} \ndef{\frQ}{{\mathfrak Q}}
\ndef{\clR}{{\mathcal R}} \ndef{\rmR}{{\mathrm R}} \ndef{\mbR}{{\mathbb R}} \ndef{\bfR}{{\mathbf R}} \ndef{\euR}{{\EuScript R}} \ndef{\frR}{{\mathfrak R}}
\ndef{\clS}{{\mathcal S}} \ndef{\rmS}{{\mathrm S}} \ndef{\mbS}{{\mathbb S}} \ndef{\bfS}{{\mathbf S}} \ndef{\euS}{{\EuScript S}} \ndef{\frS}{{\mathfrak S}}
\ndef{\clT}{{\mathcal T}} \ndef{\rmT}{{\mathrm T}} \ndef{\mbT}{{\mathbb T}} \ndef{\bfT}{{\mathbf T}} \ndef{\euT}{{\EuScript T}} \ndef{\frT}{{\mathfrak T}}
\ndef{\clU}{{\mathcal U}} \ndef{\rmU}{{\mathrm U}} \ndef{\mbU}{{\mathbb U}} \ndef{\bfU}{{\mathbf U}} \ndef{\euU}{{\EuScript U}} \ndef{\frU}{{\mathfrak U}}
\ndef{\clV}{{\mathcal V}} \ndef{\rmV}{{\mathrm V}} \ndef{\mbV}{{\mathbb V}} \ndef{\bfV}{{\mathbf V}} \ndef{\euV}{{\EuScript V}} \ndef{\frV}{{\mathfrak V}}
\ndef{\clW}{{\mathcal W}} \ndef{\rmW}{{\mathrm W}} \ndef{\mbW}{{\mathbb W}} \ndef{\bfW}{{\mathbf W}} \ndef{\euW}{{\EuScript W}} \ndef{\frW}{{\mathfrak W}}
\ndef{\clX}{{\mathcal X}} \ndef{\rmX}{{\mathrm X}} \ndef{\mbX}{{\mathbb X}} \ndef{\bfX}{{\mathbf X}} \ndef{\euX}{{\EuScript X}} \ndef{\frX}{{\mathfrak X}}
\ndef{\clY}{{\mathcal Y}} \ndef{\rmY}{{\mathrm Y}} \ndef{\mbY}{{\mathbb Y}} \ndef{\bfY}{{\mathbf Y}} \ndef{\euY}{{\EuScript Y}} \ndef{\frY}{{\mathfrak Y}}
\ndef{\clZ}{{\mathcal Z}} \ndef{\rmZ}{{\mathrm Z}} \ndef{\mbZ}{{\mathbb Z}} \ndef{\bfZ}{{\mathbf Z}} \ndef{\euZ}{{\EuScript Z}} \ndef{\frZ}{{\mathfrak Z}}

\ndef{\tA}{{\widetilde A}} \ndef{\tcA}{{\widetilde\clA}} \ndef{\ttcA}{\widetilde{\tcA}} \ndef{\sfA}{{\textsf A}} \ndef{\ttA}{\widetilde{\tA}} \ndef{\dzA}{{A^\sharp}}
\ndef{\tB}{{\widetilde B}} \ndef{\tcB}{{\widetilde\clB}} \ndef{\ttcB}{\widetilde{\tcB}} \ndef{\sfB}{{\textsf B}} \ndef{\ttB}{\widetilde{\tB}} \ndef{\dzB}{{B^\sharp}}
\ndef{\tC}{{\widetilde C}} \ndef{\tcC}{{\widetilde\clC}} \ndef{\ttcC}{\widetilde{\tcC}} \ndef{\sfC}{{\textsf C}} \ndef{\ttC}{\widetilde{\tC}} \ndef{\dzC}{{C^\sharp}}
\ndef{\tD}{{\widetilde D}} \ndef{\tcD}{{\widetilde\clD}} \ndef{\ttcD}{\widetilde{\tcD}} \ndef{\sfD}{{\textsf D}} \ndef{\ttD}{\widetilde{\tD}} \ndef{\dzD}{{D^\sharp}}
\ndef{\tE}{{\widetilde E}} \ndef{\tcE}{{\widetilde\clE}} \ndef{\ttcE}{\widetilde{\tcE}} \ndef{\sfE}{{\textsf E}} \ndef{\ttE}{\widetilde{\tE}} \ndef{\dzE}{{E^\sharp}}
\ndef{\tF}{{\widetilde F}} \ndef{\tcF}{{\widetilde\clF}} \ndef{\ttcF}{\widetilde{\tcF}} \ndef{\sfF}{{\textsf F}} \ndef{\ttF}{\widetilde{\tF}} \ndef{\dzF}{{F^\sharp}}
\ndef{\tG}{{\widetilde G}} \ndef{\tcG}{{\widetilde\clG}} \ndef{\ttcG}{\widetilde{\tcG}} \ndef{\sfG}{{\textsf G}} \ndef{\ttG}{\widetilde{\tG}} \ndef{\dzG}{{G^\sharp}}
\ndef{\tH}{{\widetilde H}} \ndef{\tcH}{{\widetilde\clH}} \ndef{\ttcH}{\widetilde{\tcH}} \ndef{\sfH}{{\textsf H}} \ndef{\ttH}{\widetilde{\tH}} \ndef{\dzH}{{H^\sharp}}
\ndef{\tI}{{\widetilde I}} \ndef{\tcI}{{\widetilde\clI}} \ndef{\ttcI}{\widetilde{\tcI}} \ndef{\sfI}{{\textsf I}} \ndef{\ttI}{\widetilde{\tI}} \ndef{\dzI}{{I^\sharp}}
\ndef{\tJ}{{\widetilde J}} \ndef{\tcJ}{{\widetilde\clJ}} \ndef{\ttcJ}{\widetilde{\tcJ}} \ndef{\sfJ}{{\textsf J}} \ndef{\ttJ}{\widetilde{\tJ}} \ndef{\dzJ}{{J^\sharp}}
\ndef{\tK}{{\widetilde K}} \ndef{\tcK}{{\widetilde\clK}} \ndef{\ttcK}{\widetilde{\tcK}} \ndef{\sfK}{{\textsf K}} \ndef{\ttK}{\widetilde{\tK}} \ndef{\dzK}{{K^\sharp}}
\ndef{\tL}{{\widetilde L}} \ndef{\tcL}{{\widetilde\clL}} \ndef{\ttcL}{\widetilde{\tcL}} \ndef{\sfL}{{\textsf L}} \ndef{\ttL}{\widetilde{\tL}} \ndef{\dzL}{{L^\sharp}}
\ndef{\tM}{{\widetilde M}} \ndef{\tcM}{{\widetilde\clM}} \ndef{\ttcM}{\widetilde{\tcM}} \ndef{\sfM}{{\textsf M}} \ndef{\ttM}{\widetilde{\tM}} \ndef{\dzM}{{M^\sharp}}
\ndef{\tN}{{\widetilde N}} \ndef{\tcN}{{\widetilde\clN}} \ndef{\ttcN}{\widetilde{\tcN}} \ndef{\sfN}{{\textsf N}} \ndef{\ttN}{\widetilde{\tN}} \ndef{\dzN}{{N^\sharp}}
\ndef{\tO}{{\widetilde O}} \ndef{\tcO}{{\widetilde\clO}} \ndef{\ttcO}{\widetilde{\tcO}} \ndef{\sfO}{{\textsf O}} \ndef{\ttO}{\widetilde{\tO}} \ndef{\dzO}{{O^\sharp}}
\ndef{\tP}{{\widetilde P}} \ndef{\tcP}{{\widetilde\clP}} \ndef{\ttcP}{\widetilde{\tcP}} \ndef{\sfP}{{\textsf P}} \ndef{\ttP}{\widetilde{\tP}} \ndef{\dzP}{{P^\sharp}}
\ndef{\tQ}{{\widetilde Q}} \ndef{\tcQ}{{\widetilde\clQ}} \ndef{\ttcQ}{\widetilde{\tcQ}} \ndef{\sfQ}{{\textsf Q}} \ndef{\ttQ}{\widetilde{\tQ}} \ndef{\dzQ}{{Q^\sharp}}
\ndef{\tR}{{\widetilde R}} \ndef{\tcR}{{\widetilde\clR}} \ndef{\ttcR}{\widetilde{\tcR}} \ndef{\sfR}{{\textsf R}} \ndef{\ttR}{\widetilde{\tR}} \ndef{\dzR}{{R^\sharp}}
\ndef{\tS}{{\widetilde S}} \ndef{\tcS}{{\widetilde\clS}} \ndef{\ttcS}{\widetilde{\tcS}} \ndef{\sfS}{{\textsf S}} \ndef{\ttS}{\widetilde{\tS}} \ndef{\dzS}{{S^\sharp}}
\ndef{\tT}{{\widetilde T}} \ndef{\tcT}{{\widetilde\clT}} \ndef{\ttcT}{\widetilde{\tcT}} \ndef{\sfT}{{\textsf T}} \ndef{\ttT}{\widetilde{\tT}} \ndef{\dzT}{{T^\sharp}}
\ndef{\tU}{{\widetilde U}} \ndef{\tcU}{{\widetilde\clU}} \ndef{\ttcU}{\widetilde{\tcU}} \ndef{\sfU}{{\textsf U}} \ndef{\ttU}{\widetilde{\tU}} \ndef{\dzU}{{U^\sharp}}
\ndef{\tV}{{\widetilde V}} \ndef{\tcV}{{\widetilde\clV}} \ndef{\ttcV}{\widetilde{\tcV}} \ndef{\sfV}{{\textsf V}} \ndef{\ttV}{\widetilde{\tV}} \ndef{\dzV}{{V^\sharp}}
\ndef{\tW}{{\widetilde W}} \ndef{\tcW}{{\widetilde\clW}} \ndef{\ttcW}{\widetilde{\tcW}} \ndef{\sfW}{{\textsf W}} \ndef{\ttW}{\widetilde{\tW}} \ndef{\dzW}{{W^\sharp}}
\ndef{\tX}{{\widetilde X}} \ndef{\tcX}{{\widetilde\clX}} \ndef{\ttcX}{\widetilde{\tcX}} \ndef{\sfX}{{\textsf X}} \ndef{\ttX}{\widetilde{\tX}} \ndef{\dzX}{{X^\sharp}}
\ndef{\tY}{{\widetilde Y}} \ndef{\tcY}{{\widetilde\clY}} \ndef{\ttcY}{\widetilde{\tcY}} \ndef{\sfY}{{\textsf Y}} \ndef{\ttY}{\widetilde{\tY}} \ndef{\dzY}{{Y^\sharp}}
\ndef{\tZ}{{\widetilde Z}} \ndef{\tcZ}{{\widetilde\clZ}} \ndef{\ttcZ}{\widetilde{\tcZ}} \ndef{\sfZ}{{\textsf Z}} \ndef{\ttZ}{\widetilde{\tZ}} \ndef{\dzZ}{{Z^\sharp}}

\ndef{\bfc}{{\bf c}}

  \ndef{\eps}{\varepsilon}
\let\geq\geqslant
\let\leq\leqslant

\ndef{\lims}[1]{\lim\limits_{#1}}
\ndef{\sums}[1]{\sum\limits_{#1}}
\ndef{\ints}[1]{\int\limits_{#1}}
\ndef{\sups}[1]{\sup\limits_{#1}}
\ndef{\liminfty}[1]{\lims{#1\to\infty}}
\ndef{\suminf}[1]{\sums{#1=1}^\infty}

\ndef{\limo}[1]{\omega\mbox{-}\!\!\!\lims{#1\to\infty}}          % \omega limit
\ndef{\limL}[1]{\rmL\mbox{-}\!\!\!\lims{#1\to\infty}}            % "L" limit
\ndef{\limLOne}[1]{\clL_1\mbox{-}\!\!\!\lims{#1}}
\ndef{\tildelimo}[1]{\tilde\omega\mbox{-}\!\!\!\lims{#1\to\infty}}
\ndef{\slim}{\mathrm{s}\mbox{-}\!\!\lim}          % strong limit
\ndef{\wlim}{\mathrm{w}\mbox{-}\!\lim}          % strong limit

\ndef{\Aut}{\operatorname{Aut}}      % group of automorphisms
\ndef{\Ch}{\operatorname{ch}}        % Chern character
\ndef{\End}{\operatorname{End}}      % group of endomorphisms
\ndef{\Hom}{\operatorname{Hom}}      % group of homomorphisms
\ndef{\Ker}{\operatorname{Ker}}      % kernel of operator
\ndef{\Log}{\operatorname{Log}}      % logarithm
\ndef{\OP}{\operatorname{OP}}        % abstract PDO's
\ndef{\Op}{\operatorname{Op}}        % abstract PDO's
\ndef{\Symb}{\operatorname{Symb}}    % symbol
\ndef{\Tr}{\operatorname{Tr}}        % usual trace
\ndef{\Wres}{\operatorname{Wres}}    % Wodzicki residue
\ndef{\cl}{\operatorname{cl}}        % Clifford
\ndef{\com}{\operatorname{com}}
\ndef{\const}{\operatorname{const}}  % constant
\ndef{\conv}{\operatorname{conv}}    % convex hull
\rndef{\det}{\operatorname{det}}     % determinant
\ndef{\detFK}[1]{\Delta\brs{#1}} % Fuglede-Kadison's determinant
\ndef{\detFKrel}[2]{\Delta_{#2}\brs{#1}} % Fuglede-Kadison's determinant

\ndef{\diag}{\operatorname{diag}}    % diagonal operator
\ndef{\dist}{\operatorname{dist}}    % distance
\ndef{\dom}{\operatorname{dom}}      % domain
\ndef{\ec}{\operatorname{ec}}        % essential codimension
\ndef{\id}{1}                        % identity operator
\ndef{\ind}{\operatorname{ind}}      % index
\ndef{\mydeg}{\operatorname{deg}}    % degree of a diff. form
\ndef{\op}{\operatorname{op}}
\ndef{\rank}{\operatorname{rank}}
\ndef{\res}{\operatorname{res}}      % residue
\ndef{\rng}{\operatorname{ran}}      % range
\ndef{\sflow}{\operatorname{sf}}     % spectral flow
\ndef{\isf}{\operatorname{isf}}      % infinitesimal spectral flow
\ndef{\sign}{\operatorname{sign}}    % signum (a la C.-Ph.)
\ndef{\sgn}{\operatorname{sgn}}      % signum (a la Connes)
\ndef{\sing}{\operatorname{sing}}    % singular
\ndef{\supp}{\operatorname{supp}}    % support
\ndef{\tr}{\operatorname{tr}}        % trace
\ndef{\var}{\operatorname{var}}      % variation of measure
\ndef{\vol}{\operatorname{vol}}      % volume or volume form
\ndef{\wn}{\operatorname{wn}}        % winding number
\ndef{\wres}{\operatorname{wres}}    % Wodzicki residue
\rndef{\Im}{\operatorname{Im}}       % imaginary part of an operator
\rndef{\Re}{\operatorname{Re}}       % real part of an operator

\ndef{\prng}[1]{\mathrm R_{#1}} % {[\rng {#1}]}          % projection onto the range of #1
\ndef{\pker}[1]{\mathrm N_{#1}} % {[\ker {#1}]}          % projection onto the kernel of #1
\ndef{\rprng}[2]{\mathrm R_{#1}^{#2}}           % projection onto the relative range of #1
\ndef{\rpker}[2]{\mathrm N_{#1}^{#2}}           % projection onto the relative kernel of #1
\ndef{\rsupp}[1]{\supp_r(#1)}
\ndef{\lsupp}[1]{\supp_l(#1)}
\ndef{\rslv}[1]{R_z(#1)}      % resolvent
\ndef{\HH}{H}                 % initial operator of perturbation theory
\ndef{\tHH}{\tilde \HH}       % final operator of perturbation theory
\ndef{\VV}{V}                 % perturbation operator of perturbation theory
\ndef{\Rz}{R_z}               % resolvent of the initial operator
\ndef{\tRz}{\tR_z}            % resolvent of the final operator
\ndef{\psif}[1]{#1^{[1]}} % {\psi_{#1}}  % divided difference
\ndef{\CPlus}[1]{W_{#1}(\mbR)}
\ndef{\bndl}{\xi}                         % vector bundle
\ndef{\bndlA}{\eta}                       % vector bundle
\ndef{\GlueMap}{\varphi}                  % glue map of a bundle
\ndef{\ChartMap}{h}                       % chart diffeomorphism map of a manifold

\ndef{\hilb}{\clH}                     % Hilbert space
   \ndef{\hilbasargument}{(\hilb)} %{(\hilb)}
\ndef{\LpH}[1]{\clL^{#1}\hilbasargument}       % the set of ...
\ndef{\saLpH}[1]{\clL_{sa}^{#1}\hilbasargument}       % the set of ...
\ndef{\clBH}{\clB\hilbasargument}              % the set of BOUNDED linear operators on Hilbert space
\ndef{\ubBH}{\clB_1\hilbasargument}            % the unit ball of the algebra of all BOUNDED linear operators on Hilbert space
\ndef{\clCH}{\clC\hilbasargument}              % the set of CLOSED DENSELY-DEFINED linear operators on Hilbert space
\ndef{\clKH}{\clK\hilbasargument}              % the set of COMPACT operators
\ndef{\clFH}{\clF\hilbasargument}              % the set of BOUNDED FREDHOLM operators
\ndef{\clUH}{\clU\hilbasargument}              % the set of UNITARIES on Hilbert space
\ndef{\clCFH}{{\clC\clF}\hilbasargument}       % the set of CLOSED DENSELY-DEFINED FREDHOLM OPERATORS on Hilbert space
\ndef{\saBH}{\clB_{sa}\hilbasargument}         % the set of S.-A. BOUNDED operators on Hilbert space
\ndef{\saCH}{\clC_{sa}\hilbasargument}         % the set of CLOSED DENSELY-DEFINED S.-A. operators on Hilbert space
\ndef{\saFH}{\clF_{sa}\hilbasargument}         % the set of BOUNDED FREDHOLM s.-a. operators
\ndef{\saKH}{\clK_{sa}\hilbasargument}         % the set of COMPACT S.-A. operators
\ndef{\saCFH}{\clC\clF_{sa}\hilbasargument}    % the set of CLOSED DENSELY-DEFINED S.-A. FREDHOLM operators on Hilbert space
\ndef{\clUFH}{\clU\clF\hilbasargument}         % the set of UNITARIES such that U+I is FREDHOLM
\ndef{\Uinj}{\clU_{inj}\hilbasargument}        % the set of UNITARIES such that U-I is injective
\ndef{\UFinj}{\clU\clF_{inj}\hilbasargument}   % the set of UNITARIES such that U-I is INJECTIVE and U+I is FREDHOLM

\ndef{\spproj}[2]{E^{#1}_{#2}}                      % spectral projection of #1
\ndef{\spprojb}[2]{E^{#2}_{#1}}                     % spectral projection of #1

\ndef{\LpN}[1]{\clL^{#1}(\clN,\tau)}     % noncommutative \mathcal L_p space
\ndef{\saLpN}[1]{\clL^{#1}_{sa}(\clN,\tau)} % s.-a. part of noncommutative \mathcal L_p space
\ndef{\rLpN}[1]{L^{#1}(\clN,\tau)}       % noncommutative L_p space
\ndef{\clAND}{(\clA,\clN,D)}             % spectral triple (A,N,D)
\ndef{\clBA}{{\clB(\clA)}}
\ndef{\saKN}{{\clK_{sa}(\clN,\tau)}}          % s.-a. \tau-compact operators
\ndef{\clKN}{{\clK(\clN,\tau)}}          % \tau-compact operators
\ndef{\clKtN}{{\clK(\tilde\clN,\tau)}}   % \tau-compact (maybe unbounded) operators
\ndef{\clFN}{{\clF(\clN,\tau)}}          % \tau-Fredholm operators
\ndef{\saFN}{{\clF_{sa}(\clN,\tau)}}     % self-adjoint \tau-Fredholm operators
\ndef{\clPN}{\clP(\clN)}                 % projections of N
\ndef{\clQN}{\clQ(\clN,\tau)}            % Calkin algebra N/K
\ndef{\infPN}{{\clP_\tau^\infty(\clN)}}  % infinite projections of N
\ndef{\clOF}[2]{\clF_{#1\mbox{-}#2}(\clN,\tau)}         % relatively Fredholm operators
\ndef{\oind}[2]{{\rm \tau\mbox{-}ind}_{#1\mbox{-}#2}}   % relative index
\ndef{\tind}{\tau\mbox{-}\ind}                  % semifinite index
\ndef{\DInd}{\ind_{\clD,\tau}}           % semifinite index
\ndef{\BF}{Breuer-Fredholm}              % Breuer-Fredholm
\ndef{\skewfred}[2]{$(#1\cdot #2)$ $\tau$\tire Fredholm}   % skew corner Fredholm
\ndef{\affl}{\eta}                       % affiliated
\ndef{\vNa}{von Neumann algebra}         % von Neumann algebra
\ndef{\nsf}{faithful normal semifinite } % normal semifinite faithful
\ndef{\taubrs}[1]{\tau\brackets{#1}}     % n.s.f. trace with brackets
\ndef{\sqbrs}[1]{\left[#1\right]}        % brackets
\ndef{\domd}{\bigcap\limits_{n\ge 0} \dom\;\delta^n}         % domain of \delta^n's
\ndef{\DiffOP}{{\rm \clD}}
\ndef{\ADA}{\clA \cup [\clD,\clA]}
\ndef{\DixIdeal}[1]{\LpH{#1,\infty}}               % Dixmier ideal
\ndef{\dixideal}{\ell^{1,\infty}}                  % commutative Dixmier ideal
\ndef{\WDixIdeal}{\LpH{1,\mathrm w}}               % weak Dixmier ideal
\ndef{\DixIdealPos}[1]{\DixIdeal{#1}_+}            % positive part of Dixmier ideal
\ndef{\DixIdealN}[1]{\LpN{#1,\infty}}              % semifinite Dixmier ideal
\ndef{\DixIdealNPar}[2]{\clL^{#1,\infty}_{#2}(\clN,\tau)}    % semifinite Dixmier ideal
\ndef{\DixIdealNPos}[1]{\LpN{#1,\infty}_+}                   % positive part of semifinite Dixmier ideal
\ndef{\TrD}{\Tr_\omega}                                      % Dixmier trace
\ndef{\tauD}{{\tau_\omega}}                                  % semifinite Dixmier trace
\ndef{\ILog}{\frac 1{\log(1+t)}}
\ndef{\ILogN}{\frac 1{\log(1+N)}}
\ndef{\DixNorm}[1]{\norm{#1}_{(1,\infty)}}                   % Dixmier norm
\ndef{\DixInt}[1]{\ints 0^t \mu_s(#1)\,ds}
\ndef{\DixIntL}[1]{\ints 0^{\lambda_{1/t}(#1)}\mu_s(#1)\,ds}
    \ndef{\SmallIdeal}{{\clL^{1, \mathrm w}}}
    \ndef{\SmallIdealMeas}{{\clL^{1, \mathrm w}_m}}
    \ndef{\DixIntII}[1]{\ints 0^t \mu_s(#1)\,ds}
    \ndef{\DixIntf}[1]{f_t(#1)}
    \ndef{\DixIntg}[1]{g_t(#1)}

\ndef{\lpi}{\clL^{1,\pi}(\clN,\tau)}
\ndef{\IIinfty}{$\mathrm{II}_\infty$\ }

\ndef{\fourier}[1]{\clF(#1)}          % Fourier transform of #1
\ndef{\HaarMeasBohrs}{\nu}            % Haar measure of the Bohr compact
\ndef{\BrownsMeas}{\mu}               % Brown's measure
\ndef{\BohrCont}[1]{\tilde{#1}}       % continuation of a function to the Bohr compact
\ndef{\APMean}{{M}}                   % mean value of a.p. function
\ndef{\CDSS}{{\clA_B}}                % Coburn-Douglas-Schaeffer-Singer's factor
\ndef{\matr}{{\rm Mat}}               % standard matrix algebra
\ndef{\seque}[1]{\ensuremath{\{#1_j\}_{j=1}^\infty}}    % sequence of numbers  a_1, a_2, ...
\ndef{\sequen}[2]{\ensuremath{\{#1_#2\}_{#2=1}^\infty}}    % sequence of numbers  a_1, a_2, ...
\ndef{\Seque}[1]{\ensuremath{\left(#1_0,#1_1,#1_2,\dots\right)}}    % sequence of numbers  a_1, a_2, ...
\ndef{\Cesaro}{H}                           % the Cesaro operator (on sequences)
\ndef{\CesaroRPlus}{M}                      % the Cesaro operator on positive semiaxis
\ndef{\Dilation}{D}                         % the dilation operator (on sequences)
\ndef{\Shift}{T}                            % the shift operator (on sequences)

\ndef{\norm}[1]{\left\Vert#1\right\Vert}    % norm of #1
\ndef{\TrNorm}[1]{\norm{#1}_1}              % trace norm of #1
\ndef{\HSNorm}[1]{\norm{#1}_2}              % Hilbert-Schmidt norm of #1
\ndef{\InftyNorm}[1]{\norm{#1}_\infty}      % uniform norm of #1
\ndef{\normQN}[1]{\norm{#1}_{\clQN}}        % Calkin norm of #1
\ndef{\clLnorm}[1]{\norm{#1}_{1,\infty}}    % $1,\infty$- trace norm of #1

\ndef{\ccurve}{\gamma}                      % a curve in complex plane for Cauchy integral

\ndef{\abs}[1]{\left\lvert#1\right\rvert}   % absolute value of #1
\ndef{\set}[1]{\left\{#1\right\}}           % set of ...
\ndef{\brackets}[1]{\left(#1\right)}        % brackets
\ndef{\brs}[1]{\brackets{#1}}               % brackets
\ndef{\Brs}[1]{\big(#1\big)}                % brackets
\ndef{\BRS}[1]{\Big(#1\Big)}                % brackets
\ndef{\scal}[2]{\left\la #1,#2\right\ra}               % scalar product
\ndef{\precprec}{\prec\!\!\!\prec}
\ndef{\qeq}{\stackrel?=}
\ndef{\spectrum}[1]{\sigma_{#1}} %{\mathrm{Spec}(#1)}       % spectrum of an operator
\ndef{\numrange}[1]{\mathrm{W}(#1)}                         % numerical range of an operator
\rndef{\emptyset}{\varnothing}                              % empty set
\ndef{\csupp}{c}                           % subscript for compactly supported functions
\ndef{\closure}[1]{\overline{#1}}
\ndef{\linspan}[1]{\mathrm{span}\ {#1}}
\ndef{\bddborel}[1]{B(#1)}                 % the space of bounded Borel functions on the measure space #1
\ndef{\charfunc}{\chi}
\rndef{\ln}{\log}
\ndef{\FrDer}{\euD}                        % Fr\'echet derivative
\ndef{\LieDer}[1]{\pounds_{#1}\,}          % Lie derivative
\ndef{\dds}{\left.\frac d{ds} \right|_{s = 0}}
\ndef{\ortcmp}[1]{#1^{\scriptscriptstyle \perp}}            % orthogonal complement of projection #1
\ndef{\Laplace}{\Delta}                    % Laplace operator

\ndef{\matrPQ}[3]
{
    \left(
      \begin{array}{cc}
        #1_{11} & #1_{12} \\
        #1_{21} & #1_{22}
      \end{array}
    \right)_{[#2,#3]}
}
\newcounter{margcomcount}
\ndef{\margcom}[1]{\marginpar{\bf \small #1} \addtocounter{margcomcount}{1}}

\newcounter{margproof}
\ndef{\margproof}{\marginpar{\bf \small Proof} \addtocounter{margproof}{1}}

\newcounter{margdetails}
\ndef{\margdetails}{\marginpar{\bf \qquad\qquad \small details} \addtocounter{margdetails}{1}}

\newcounter{margproofb}
\ndef{\margproofb}{\marginpar{\bf \small Proof (B)} \addtocounter{margproof}{1}}

\newcounter{margdetailsb}
\ndef{\margdetailsb}{\marginpar{\bf \small Details (B)} \addtocounter{margdetailsb}{1}}

\ndef{\mytimes}{\!\times\!}
\ndef{\sss}[1]{\subsubsection{}\label{#1}}
\rndef{\phi}{\varphi}
\ndef{\OpenUnitDisk}{D}
\ndef{\RHS}{RHS}                            % right hand side
\ndef{\LHS}{LHS} %right and left hand side  % left hand side
\ndef{\ttt}{\Leftrightarrow}
\ndef{\then}{\Rightarrow}
\ndef{\tto}{\longrightarrow}
\ndef{\nno}{\nonumber\\}
\ndef{\newn}[1]{\index{#1} \emph{#1}}       % new notion
\ndef{\la}{\langle}
\ndef{\ra}{\rangle}
\ndef{\dbar}{{\;\bar{\phantom{o}} \!\!\!\! d}}
\ndef{\stl}[1]{\stackrel{\vbox to 0pt{\vss\hbox{$\scriptstyle #1$}}}}
\ndef{\mathcomment}[1]{{\scriptstyle\text{(#1)}}\qquad}        % for comments at the ends of lines with math formulas
\ndef{\details}[1]{\smallskip\begin{center} {\bf Here:} #1\end{center}\medskip}
\ndef{\indexcom}[1]{ --- #1}
\ndef{\longsim}{\ \sim \ }              % for use in formulas.
\ndef{\tire}{-}              % for use in formulas.
\ndef{\intinfinf}{\int_{-\infty}^\infty}
\ndef{\refnsftrace}{\cite[V.\,2.\,1]{TakI}} % reference to definition of n.f.s. trace %\cite[Definition V.2.1]{TakI}
\ndef{\refaffloper}{\cite[IV.\,5, Exercise 3]{TakI}} % reference to definition of affiliated operator
\ndef{\refsemifinvNa}{\cite[V.\,1.\,21]{TakI}} %  semifinite vNa
\ndef{\reftaumeasurable}{\cite[Definition 1.2]{FK86PJM}} % tau-measurable operator
\ndef{\reftautraceclassaffl}{\cite[V.2, p.\,320]{TakI}} % tau-trace class affiliated operator
\ndef{\refinvoperideal}{\cite[Appendix A.2]{CP2}} % invariant operator ideal
\ndef{\reftautracenorm}{\cite[V.2, p.\,320]{TakI}} % tau trace norm (1-norm)
\ndef{\reftaucompact}{\cite{}} % tau compact operator
\ndef{\reftauFredholm}{\cite[Appendix B]{PR94JFA}} % tau Fredholm operator

     \ndef{\npartial}{\slash\!\!\!\partial}
     \ndef{\Heis}{\operatorname{Heis}}
     \ndef{\Solv}{\operatorname{Solv}}
     \ndef{\Spin}{\operatorname{Spin}}
     \ndef{\SO}{\operatorname{SO}}
     \ndef{\Index}{\operatorname{index}}
             \ndef{\coker}{{\mbox coker}}
             \ndef{\p}{\partial}
             \ndef{\dd}{|\clD|}
             \ndef{\n}{\parallel}

     \ndef{\gf}[2]{\genfrac{}{}{0pt}{}{#1}{#2}}
     \ndef{\ta}{\widetilde{\alpha}}
     \ndef{\tb}{\widetilde{\beta}}
     \ndef{\txi}{\widetilde{\xi}}
     \ndef{\tk}{\widetilde{K}}
     \ndef{\CGh}{\widetilde{\CG}}
     \ndef{\boe}{{\bf e}}\ndef{\bt}{{\bf t}}
     \ndef{\vth}{\vartheta}
     \ndef{\db}{\overline{\partial}}
     \ndef{\hV}{\hat{V}}
     \ndef{\cag}{{\clA^\Gamma}}
     \ndef{\sind}{\sigma{\rm -ind}}

\let\LatexCite=\cite  % just renaming

\let\ifnumref\iftrue % use this command if you want NUMBER REFERENCES like [4]
\ndef{\ifuncited}[4]{\expandafter\ifx\csname used#4\endcsname\relax}

\ndef{\ifcited}[4]{\expandafter\ifx\csname used#4\endcsname\relax\else}
  \ndef{\papertitle}[1]{ \emph{#1}, }
  \ndef{\paperauthor}[2]{#2}  %{\ifnum#1=0$^*$\fi#2}
  \ndef{\pbbi}[9]{%
      \ifcited{#1}{#2}{#3}{#5}%
        \ifnumref%
          \bibitem{#5}\paperauthor{#1}{#6},\papertitle{#7}#8.%
        \else%
          \advance #9 by 1%
          \ifnum#9<1%
            \bibitem[#4]{#5}\paperauthor{#1}{#6}, \papertitle{#7}#8.%
          \else%
            \bibitem[#4$\!_{\the#9}\!$]{#5}\paperauthor{#1}{#6},\papertitle{#7}#8.%
          \fi%
        \fi%
      \fi%
  }
  \ndef{\mbbi}[8]{%
     \ifcited{#1}{#2}{#3}{#5}%
        \ifnumref%
          \bibitem{#5}\paperauthor{#1}{#6},\papertitle{#7}#8.%
        \else%
          \bibitem[#4]{#5}\paperauthor{#1}{#6},\papertitle{#7}#8.%
        \fi%
     \fi%
  }
\ndef{\AddCite}[1]{%
   \ifuncited{0}{0}{0}{#1}%
     \expandafter\gdef\csname used#1\endcsname {}%
   \fi%
}

\def\ProcessCite#1,{%
     \ifx\relax#1%
         \let\next=\relax%
     \else%
         \AddCite{#1}%
         \let\next=\ProcessCite%
     \fi%
     \next%
}

\ndef{\AddCites}[1]{\ProcessCite#1,\relax,}

\ndef{\CiteWithoutExtension}[1]{%
   \AddCites{#1}%
   \LatexCite{#1}%
}

\def\CiteWithExtension[#1]#2{%
   \AddCites{#2}%
   \LatexCite[#1]{#2}%
}

\ndef{\CleverCite}{%
    \ifx\NChar[ %
       \let\MyCite=\CiteWithExtension %
    \else %
       \let\MyCite=\CiteWithoutExtension %
    \fi %
    \MyCite%
}

\renewcommand{\cite}{\futurelet\NChar\CleverCite}
      \ndef{\volume}[1]{{\bf #1}}
      \ndef{\VolYearPP}[3]{\ifnum#2=0 (to appear)\else\volume{#1} (#2), #3\fi}
      \ndef{\VolNoYearPP}[4]{\ifnum#3=0 (to appear)\else\volume{#1} #2 (#3), #4\fi}
      \ndef{\libcode}[1]{}%{{,\bf\ #1}}

\ndef{\jnActaMath}[3]{Acta Math. \VolYearPP{#1}{#2}{#3}}                       % ActM        \libcode{510.5 A18}
\ndef{\jnAdvMath}[3]{Adv. in~Math. \VolYearPP{#1}{#2}{#3}}                     % AdvM        \libcode{510.5 A24}
\ndef{\jnAlgAnal}[3]{Algebra i~Analiz \VolYearPP{#1}{#2}{#3}}
\ndef{\jnAmerMathMonth}[3]{Amer. Math. Monthly \VolYearPP{#1}{#2}{#3}}         % AMM         \libcode{510.5 A3}
\ndef{\jnAnnMath}[4]{Ann. of~Math. \VolNoYearPP{#1}{#2}{#3}{#4}}               % AnnM        \libcode{510.5 A61}
\ndef{\jnAnalMath}[3]{J. Anal. Math. \VolYearPP{#1}{#2}{#3}}                   % JAnalM      \libcode{}
\ndef{\jnBullLondMathSoc}[3]{Bull. London Math. Soc. \VolYearPP{#1}{#2}{#3}}   % BLMS        \libcode{}
\ndef{\jnBullAMS}[3]{Bull. Amer. Math. Soc. \VolYearPP{#1}{#2}{#3}}   % BLMS        \libcode{}
\ndef{\jnCanMathBull}[3]{Canad. Math. Bull. \VolYearPP{#1}{#2}{#3}}            % CMB         \libcode{}
\ndef{\jnCanMath}[3]{Canad. J.~Math. \VolYearPP{#1}{#2}{#3}}             % CJM         \libcode{}
\ndef{\jnCommMathPhys}[3]{Comm. Math. Phys \VolYearPP{#1}{#2}{#3}}             % CMP         \libcode{530.1505 C73}
\ndef{\jnCommPDE}[3]{Comm. Partial Differential Equations \VolYearPP{#1}{#2}{#3}}             % CMP         \libcode{530.1505 C73}
\ndef{\jnComptRendue}[3]{C.\,R.~Acad. Sci. Paris S\'er. A-B \VolYearPP{#1}{#2}{#3}}      % CR   \libcode{505 A 161}
\ndef{\jnDiffGeom}[3]{J.~Diff. Geom. \VolYearPP{#1}{#2}{#3}}                   % JDG         \libcode{516.3605 J86}
\ndef{\jnErgodicTheory}[3]{Ergodic Theory and Dynamical Systems \VolYearPP{#1}{#2}{#3}} % ETDS  \libcode{}
\ndef{\jnFuncAnal}[3]{J.~Functional Analysis \VolYearPP{#1}{#2}{#3}}           % JFA         \libcode{515.05 J88}
\ndef{\jnFunkAnalPril}[4]{Функциональный анализ и его приложения \VolNoYearPP{#1}{#2}{#3}{#4}}  % JFAP
\ndef{\jnGAFA}[3]{GAFA \VolYearPP{#1}{#2}{#3}}                                 % GAFA        \libcode{}
\ndef{\jnIHES}[3]{IHES Publ. Math. (Paris) \VolYearPP{#1}{#2}{#3}}             % IHES        \libcode{}
\ndef{\jnIEOT}[3]{Integral Equations Operator Theory   \VolYearPP{#1}{#2}{#3}} %         \libcode{515.4505 I61}
\ndef{\jnIsrMath}[3]{Israel J.~Math. \VolYearPP{#1}{#2}{#3}}                   % IMJ         \libcode{}
\ndef{\jnKTheory}[3]{K-Theory \VolYearPP{#1}{#2}{#3}}                          % KT          \libcode{}
\ndef{\jnLetMathPhys}[3]{Lett. Math. Phys. \VolYearPP{#1}{#2}{#3}}             % LMP         \libcode{}
\ndef{\jnMathAnn}[3]{Math. Ann. \VolYearPP{#1}{#2}{#3}}                        % MAnn        \libcode{510.5 M51}
\ndef{\jnMathAnalAppl}[3]{J.~Math. Anal. and Appl. \VolYearPP{#1}{#2}{#3}}     % MathAnalAppl  \libcode{510.5 J88} only in repository
\ndef{\jnMathNachr}[3]{Math. Nachr. \VolYearPP{#1}{#2}{#3}}
\ndef{\jnMathPhys}[3]{J. Math. Phys. \VolYearPP{#1}{#2}{#3}}
\ndef{\jnOperTheory}[3]{J.~Operator Theory \VolYearPP{#1}{#2}{#3}}             % JOT         \libcode{}
\ndef{\jnPacJMath}[3]{Pacific J.~Math. \VolYearPP{#1}{#2}{#3}}                  % PJM         \libcode{510.5 P11}
\ndef{\jnPositivity}[3]{Positivity \VolYearPP{#1}{#2}{#3}}
\ndef{\jnProcAmerMS}[3]{Proc. Amer. Math. Soc. \VolYearPP{#1}{#2}{#3}}         % PAMS        \libcode{510.5 A5p}
\ndef{\jnProcCambPhilSoc}[3]{Math. Proc. Camb. Phil. Soc. \VolYearPP{#1}{#2}{#3}}
\ndef{\jnReineAngew}[3]{J.~Reine Angew. Math. \VolYearPP{#1}{#2}{#3}}          % JRAM        \libcode{510.5 J86}
\ndef{\jnTokyoMath}[3]{Tokyo J.~Math. \VolYearPP{#1}{#2}{#3}}
\ndef{\jnTopology}[3]{Topology \VolYearPP{#1}{#2}{#3}}
\ndef{\jnTransAmerMathSoc}[3]{Trans. Amer. Math. Soc. \VolYearPP{#1}{#2}{#3}}
\ndef{\jnIzvANSSSR}[3]{Izv. Akad. Nauk SSSR, Ser. Mat. \VolYearPP{#1}{#2}{#3}}
\ndef{\jnIzvVyshUchZav}[3]{Izv. Vyssh. Uch. Zav., Mat. \VolYearPP{#1}{#2}{#3} (Russian)}
\ndef{\jnIzdatLenUniv}[2]{Izdat. Leningrad. Univ., Leningrad, (#1), #2 (Russian)}
\ndef{\jnFieldsInsComm}[3]{Fields Inst. Comm. \VolYearPP{#1}{#2}{#3}}
\ndef{\jnDoklANSSSR}[3]{Dokl. Akad. Nauk SSSR \VolYearPP{#1}{#2}{#3}}
\ndef{\jnMatZametki}[3]{Matem. zametki \VolYearPP{#1}{#2}{#3}}
\ndef{\jnRussMathSurvey}[3]{Russian Math. Surveys \VolYearPP{#1}{#2}{#3}}
\ndef{\jnSibMathJ}[3]{Sib. Math.~J. \VolYearPP{#1}{#2}{#3}}
\ndef{\jnSovMath}[3]{J.~Soviet math. \VolYearPP{#1}{#2}{#3}}
\ndef{\jnTransMoscMathSoc}[3]{Trans. Moscow Math. Soc. \VolYearPP{#1}{#2}{#3}}
\ndef{\jnUMN}[3]{Uspekhi Mat. Nauk \VolYearPP{#1}{#2}{#3}}

\ndef{\bkTransMathMon}[2]{Trans. Math. Monographs, AMS, \volume{#1}, #2}

\ndef{\pbBirkhauser}[1]{Birkh\"auser, Boston, #1}
\ndef{\pbFactorial}[1]{Moscow, Factorial, #1}
\ndef{\pbGauthier}[1]{Gauthier-Villars, Paris, #1}
\ndef{\pbNauka}[1]{Moscow, Nauka, #1 (Russian)}
\ndef{\pbNaukaR}[1]{Москва, Наука, #1}
\ndef{\pbPrinceton}[1]{Princeton University Press, Princeton, New Jersey, #1}
\ndef{\pbPublPerish}[1]{Publish or Perish Inc., Berkeley, #1}
\ndef{\pbSpringer}[1]{Springer-Verlag, #1}

\ndef{\myauthor}[1]{\mbox{#1}}

\ndef{\Ahiezer}{\myauthor{N.\,I.\,Ahiezer }}
\ndef{\Arazy}{\myauthor{J.\,Arazy}}
\ndef{\Astashkin}{\myauthor{S.\,V.\,Astashkin}}
\ndef{\Atiyah}{\myauthor{M.\,Atiyah}}
\ndef{\Avron}{\myauthor{J.\,Avron}}
\ndef{\Azamov}{\myauthor{N.\,A.\,Azamov}}
\ndef{\Banach}{\myauthor{S.\,Banach}}
\ndef{\Benameur}{\myauthor{M-T.\,Benameur}}
\ndef{\Bennett}{\myauthor{C.\,Bennett}}
\ndef{\Berezin}{\myauthor{F.\,A.\,Berezin}}
\ndef{\Berline}{\myauthor{N.\,Berline}}
\ndef{\Birman}{\myauthor{M.\,Sh.\,Birman}}
\ndef{\Blackadar}{\myauthor{B.\,Blackadar}}
\ndef{\Bogolyubov}{\myauthor{N.\,N.\,Bogolyubov}}
\ndef{\Bonsall}{\myauthor{F.\,F.\,Bonsall}}
\ndef{\BoosBavnbek}{\myauthor{B.\,Boo$\beta$-Bavnbek}}
\ndef{\Bott}{\myauthor{R.\,Bott}}
\ndef{\Bratteli}{\myauthor{O.\,Bratteli}}
\ndef{\Bredon}{\myauthor{G.\,E.\,Bredon}}
\ndef{\Breuer}{\myauthor{M.\,Breuer}}
\ndef{\Brown}{\myauthor{L.\,G.\,Brown}}
\ndef{\Bruneau}{\myauthor{V.\,Bruneau}}
\ndef{\Buslaev}{\myauthor{V.\,S.\,Buslaev}}
\ndef{\Carey}{\myauthor{A.\,L.\,Carey}}
\ndef{\CareyRW}{\myauthor{R.\,W.\,Carey}} %Richard
\ndef{\Cartan}{\myauthor{H.\,Cartan}}
\ndef{\Chilin}{\myauthor{V.\,I.\,Chilin}}
\ndef{\Coburn}{\myauthor{L.\,A.\,Coburn}}
\ndef{\Connes}{\myauthor{A.\,Connes}}
\ndef{\Cornfeld}{\myauthor{I.\,P.\,Cornfeld}}
\ndef{\Daletskii}{\myauthor{Yu.\,L.\,Daletski\u\i}}   %Daletskii
\ndef{\Dixmier}{\myauthor{J.\,Dixmier}}
\ndef{\DoddsPG}{\myauthor{P.\,G.\,Dodds}}
\ndef{\DoddsTK}{\myauthor{T.\,K.\,Dodds}}
\ndef{\Douglas}{\myauthor{R.\,G.\,Douglas}}
\ndef{\Dubrovin}{\myauthor{B.\,A.\,Dubrovin}}
\ndef{\Dugundji}{\myauthor{J.\,Dugundji}}
\ndef{\Duncan}{\myauthor{J.\,Duncan}}
\ndef{\Dunford}{\myauthor{N.\,Dunford}}
\ndef{\Dykema}{\myauthor{K.\,J.\,Dykema}}
\ndef{\Edwards}{\myauthor{R.\,E.\,Edwards}}
\ndef{\Eilenberg}{\myauthor{S.\,Eilenberg}}
\ndef{\Fack}{\myauthor{T.\,Fack}} %Thierry
\ndef{\Faddeev}{\myauthor{L.\,D.\,Faddeev}}
\ndef{\Farber}{\myauthor{M.\,Farber}}
\ndef{\Farforovskaya}{\myauthor{Yu.\,B.\,Farforovskaya}}
\ndef{\Federer}{\myauthor{H.\,Federer}}
\ndef{\Fedosov}{\myauthor{B.\,V.\,Fedosov}}
\ndef{\Figiel}{\myauthor{T.\,Figiel}} %Tadeush?
\ndef{\Figueroa}{\myauthor{H.\,Figueroa}}
\ndef{\Fillmore}{\myauthor{P.\,A.\,Fillmore}}
\ndef{\Fomenko}{\myauthor{A.\,T.\,Fomenko}} % Anatolii Trofimovich
\ndef{\Fomin}{\myauthor{S.\,V.\,Fomin}}
\ndef{\Frohlich}{\myauthor{J.\,Fr\"ohlich}}
\ndef{\Fuglede}{\myauthor{B.\,Fuglede}}
\ndef{\Furutani}{\myauthor{K.\,Furutani}}
\ndef{\Gelfand}{\myauthor{I.\,M.\,Gelfand}}
\ndef{\Gesztesy}{\myauthor{F.\,Gesztesy}}     %Fritz
\ndef{\Getzler}{\myauthor{E.\,Getzler}} % Ezra
\ndef{\Gilkey}{\myauthor{P.\,B.\,Gilkey}}
\ndef{\Gitler}{\myauthor{S.\,Gitler}}
\ndef{\Glazman}{\myauthor{I.\,M.\,Glazman}}
\ndef{\Glimm}{\myauthor{J.\,Glimm}}
\ndef{\Gohberg}{\myauthor{I.\,C.\,Gohberg}}
\ndef{\Golze}{\myauthor{F.\,Golze}}
\ndef{\GraciaBondia}{\myauthor{J.\,M.\,Gracia-Bond\'{i}a}}
\ndef{\Greenleaf}{\myauthor{F.\,P.\,Greenleaf}}
\ndef{\Gromov}{\myauthor{M.\,Gromov}}
\ndef{\Gunning}{\myauthor{R.\,C.\,Gunning}}
\ndef{\Haagerup}{\myauthor{U.\,Haagerup}}
\ndef{\Haag}{\myauthor{R.\,Haag}}
\ndef{\Halmos}{\myauthor{Halmos}}
\ndef{\Hardy}{\myauthor{G.\,H.\,Hardy}}
\ndef{\Higson}{\myauthor{N.\,Higson}}  % Nigel
\ndef{\Hoermander}{\myauthor{L.\,Hoermander}} % Lars
\ndef{\Hoffman}{\myauthor{K.\,Hoffman}} % Kenneth Hoffman
\ndef{\Ito}{\myauthor{K.\,Ito}}
\ndef{\Jaffe}{\myauthor{A.\,Jaffe}}
\ndef{\James}{\myauthor{I.\,M.\,James}}
\ndef{\Javrjan}{\myauthor{V.\,A.\,Javrjan}}
\ndef{\Kadison}{\myauthor{R.\,V.\,Kadison}}
\ndef{\Kalton}{\myauthor{N.\,J.\,Kalton}} % Nigel
\ndef{\Kato}{\myauthor{T.\,Kato}} % Tosio
\ndef{\Kobayashi}{\myauthor{S.\,Kobayashi}}
\ndef{\Koplienko}{\myauthor{L.\,S.\,Koplienko}}
\ndef{\Korotyaev}{\myauthor{E.\,Korotyaev}}
\ndef{\Kosaki}{\myauthor{H.\,Kosaki}}
\ndef{\KreinMG}{\myauthor{M.\,G.\,Kre\u\i n}}
\ndef{\KreinSG}{\myauthor{S.\,G.\,Kre\u\i n}}
\ndef{\Leichtnam}{\myauthor{E.\,Leichtnam}}
\ndef{\Lesch}{\myauthor{M.\,Lesch}}
\ndef{\Lesniewski}{\myauthor{A.\,Lesniewski}}
\ndef{\Levitan}{\myauthor{B.\,M.\,Levitan}}
\ndef{\Lidskii}{\myauthor{V.\,B.\,Lidskii}}
\ndef{\Lifshits}{\myauthor{I.\,M.\,Lifshits}}
\ndef{\Lindenstrauss}{\myauthor{J.\,Lindenstrauss}}
\ndef{\Loday}{\myauthor{J.-L.\,Loday}}
\ndef{\Lord}{\myauthor{S.\,Lord}}      %Steven
\ndef{\Lorentz}{\myauthor{G.\,Lorentz}}
\ndef{\Magnus}{\myauthor{W.\,Magnus}}
\ndef{\Makarov}{\myauthor{K.\,A.\,Makarov}}
\ndef{\Mathai}{\myauthor{V.\,Mathai}}         %Varghese?
\ndef{\McKean}{\myauthor{H.\,P.\,McKean}}
\ndef{\Mishchenko}{\myauthor{A.\,S.\,Mishchenko}}
\ndef{\Moore}{\myauthor{C.\,C.\,Moore}}
\ndef{\Moscovici}{\myauthor{H.\,Moscovici}}  %Henri?
\ndef{\Motovilov}{\myauthor{A.\,K.\,Motovilov}}
\ndef{\Moyer}{\myauthor{R.\,D.\,Moyer}}
\ndef{\Naboko}{\myauthor{S.\,N.\,Naboko}}
\ndef{\Narasimhan}{\myauthor{R.\,Narasimhan}}
\ndef{\Nomizu}{\myauthor{K.\,Nomizu}}
\ndef{\Novikov}{\myauthor{S.\,P.\,Novikov}}
\ndef{\Osterwalder}{\myauthor{K.\,Osterwalder}}
\ndef{\Patodi}{\myauthor{V.\,Patodi}}
\ndef{\Pagter}{\myauthor{B.\,de~Pagter}}  %Ben
\ndef{\Pavlov}{\myauthor{B.\,S.\,Pavlov}}
\ndef{\Pedersen}{\myauthor{G.\,K.\,Pedersen}}
\ndef{\Peller}{\myauthor{V.\,V.\,Peller}}
\ndef{\Perera}{\myauthor{V.\,S.\,Perera}}
\ndef{\Petunin}{\myauthor{Ju.\,I.\,Petunin}}
\ndef{\Phillips}{\myauthor{J.\,Phillips}}  %John
\ndef{\Piazza}{\myauthor{P.\,Piazza}}   %Paolo
\ndef{\Pincus}{\myauthor{J.\,D.\,Pincus}}   %Joel
\ndef{\Poincare}{Poincar\'e}
\ndef{\Postnikov}{\myauthor{M.\,M.\,Postnikov}} % Mikhail
\ndef{\Prinzis}{\myauthor{R.\,Prinzis}}
\ndef{\Privalov}{\myauthor{I.\,I.\,Privalov}}
\ndef{\Pushnitski}{\myauthor{A.\,B.\,Pushnitski}} % Alexander
\ndef{\Raeburn}{\myauthor{I.\,Raeburn}}
\ndef{\Raikov}{\myauthor{G.\,Raikov}}
\ndef{\Reed}{\myauthor{M.\,Reed}}
\ndef{\Rennie}{\myauthor{A.\,Rennie}}
\ndef{\Rickart}{\myauthor{C.\,E.\,Rickart}}
\ndef{\Riesz}{\myauthor{F.\,Riesz}}
\ndef{\Ringrose}{\myauthor{J.\,Ringrose}}
\ndef{\Robinson}{\myauthor{D.\,Robinson}}
\ndef{\Rossi}{\myauthor{H.\,Rossi}}
\ndef{\Rudin}{\myauthor{W.\,Rudin}}
\ndef{\Ruelle}{\myauthor{D.\,Ruelle}}
\ndef{\Ruzhansky}{\myauthor{M.\,Ruzhansky}}
\ndef{\Sakai}{\myauthor{Sh.\,Sakai}}
\ndef{\Sargsjan}{\myauthor{I.\,S.\,Sargsjan}}
\ndef{\Sato}{\myauthor{H.\,Sato}}
\ndef{\Schaeffer}{\myauthor{D.\,G.\,Schaeffer}}
\ndef{\Schluchtermann}{\myauthor{G.\,Schluchtermann}}
\ndef{\Schochet}{\myauthor{C.\,Schochet}}
\ndef{\Schrodinger}{\myauthor{E.\,Schr\"odinger}}
\ndef{\Schrohe}{\myauthor{E.\,Schrohe}}
\ndef{\Schwartz}{\myauthor{J.\,T.\,Schwartz}}
\ndef{\Sedaev}{\myauthor{A.\,A.\,Sedaev}}
\ndef{\Seiler}{\myauthor{R.\,Seiler}}
\ndef{\Semenov}{\myauthor{E.\,M.\,Semenov}}
\ndef{\Shabat}{\myauthor{B.\,V.\,Shabat}}
\ndef{\Shafarevich}{\myauthor{I.\,R.\,Shafarevich}}
\ndef{\Sharpley}{\myauthor{R.\,Sharpley}}
\ndef{\Shilov}{\myauthor{G.\,E.\,Shilov}}
\ndef{\Shirkov}{\myauthor{D.\,V.\,Shirkov}}
\ndef{\Shubin}{\myauthor{M.\,A.\,Shubin}}
\ndef{\Silverman}{\myauthor{H.\,Silverman}}
\ndef{\Simon}{\myauthor{B.\,Simon}}
\ndef{\Sinai}{\myauthor{Ya.\,G.\,Sinai}}
\ndef{\Singer}{\myauthor{I.\,M.\,Singer}}
\ndef{\Solomyak}{\myauthor{M.\,Z.\,Solomyak}}
\ndef{\Soloviev}{\myauthor{Yu.\,P.\,Soloviev}}
\ndef{\Spivak}{\myauthor{M.\,Spivak}}
\ndef{\Stenkin}{\myauthor{V.\,V.\,Sten'kin}}
\ndef{\Stratila}{\myauthor{S.\,Stratila}}
\ndef{\Sucheston}{\myauthor{L.\,Sucheston}}
\ndef{\Sukochev}{\myauthor{F.\,A.\,Sukochev}}
\ndef{\Switzer}{\myauthor{R.\,M.\,Switzer}}
\ndef{\SzNagy}{\myauthor{B.\,Sz.-Nagy}}
\ndef{\Takesaki}{\myauthor{M.\,Takesaki}}
\ndef{\Taylor}{\myauthor{M.\,E.\,Taylor}}
\ndef{\Treves}{\myauthor{F.\,Treves}}
\ndef{\Troitsky}{\myauthor{E.\,V.\,Troitsky}}
\ndef{\Tzafriri}{\myauthor{L.\,Tzafriri}}
\ndef{\Varilly}{\myauthor{J.\,C.\,V\'{a}rilly}}
\ndef{\Vergne}{\myauthor{M.\,Vergne}}
\ndef{\Vladimirov}{\myauthor{V.\,S.\,Vladimirov}}
\ndef{\Voiculescu}{\myauthor{D.\,Voiculescu}}
\ndef{\Weiss}{\myauthor{G.\,Weiss}}
\ndef{\Wells}{\myauthor{R.\,O.\,Wells}}
\ndef{\Williams}{\myauthor{J.\,P.\,Williams}}
\ndef{\Winkler}{\myauthor{S.\,Winkler}}
\ndef{\Witten}{\myauthor{E.\,Witten}}
\ndef{\Wodzicki}{\myauthor{M.\,Wodzicki}}
\ndef{\Wojciechowski}{\myauthor{K.\,P.\,Wojciechowski}}
\ndef{\Yafaev}{\myauthor{D.\,R.\,Yafaev}}
\ndef{\Yosida}{\myauthor{K.\,Yosida}}
\ndef{\Zsido}{\myauthor{L.\,Zsido}}
\newcommand{\Cc}{C_\csupp^\infty(\mbR)}
\rndef{\HH}{H}
\rndef{\VV}{V}
\ndef{\TT}{\rmT}

\begin{document}
%\today
\title[Spectral averaging]{Spectral averaging for trace \\ compatible operators}
\author{N.\,A.\,Azamov}
\author{F.\,A.\,Sukochev}
\address{ School of Informatics and Engineering
    \\ Flinders University of South Australia \\ Bedford Park, 5042, SA Australia.
}
\email{azam0001@infoeng.flinders.edu.au, sukochev@infoeng.flinders.edu.au}
\keywords{Spectral shift function, spectral averaging, infinitesimal spectral flow,
    trace compatible operators, semifinite von Neumann algebra
}
\subjclass{ %Mathematics Subject Classification (2000).
    Primary 47A11; % Local spectral properties
    Secondary 47A55 % Perturbation theory
}

\begin{abstract}
In this note the notions of trace compatible operators
and infinitesimal spectral flow are introduced.
We define the spectral shift function as the integral of infinitesimal spectral flow.
It is proved that the spectral shift function thus defined is absolutely continuous
and Krein's formula is established.
Some examples of trace compatible affine spaces of operators are given.
\end{abstract}

\maketitle

%+++

\section{Introduction}
Let $H_0$ be a self-adjoint operator, and let $V$ be a trace class operator on a Hilbert space $\hilb.$
Then \KreinMG 's famous result \cite{Kr53MS} says that there is a unique $L^1$\tire function
$\xi_{H_0+V,H_0}(\lambda),$ known as the Krein spectral shift function, such that for any $\Cc$ function $f$
\begin{gather} \label{F: Krein's formula}
  \Tr(f(H_0+V) - f(H_0)) = \int_{-\infty}^\infty f'(\lambda) \xi_{H_0+V,H_0}(\lambda)\,d\lambda.
\end{gather}
The notion of the spectral shift function was discovered by the physicist \Lifshits\ \cite{Li52UMN}.
An excellent survey on the spectral shift function can be found in \cite{BP98IEOT}.

In 1975, Birman and Solomyak \cite{BS72SM} proved the following remarkable formula for the spectral shift function
\begin{gather} \label{F: BS formula}
  \xi(\lambda) = \frac d{d\lambda} \int_0^1 \Tr(V E_{(-\infty,\lambda]}^{H_r})\,dr,
\end{gather}
where $H_r = H_0 + rV,$ $r \in \mbR,$ and $E_{(-\infty,\lambda]}^{H_r}$ is the spectral projection.
Birman-Solomyak's proof relies on double operator integrals. An elementary derivation of
(\ref{F: BS formula}) was obtained in \cite{GMM99} (without using double operator integrals).

Actually, this spectral averaging formula was discovered for the first time by Javrjan \cite{Jav}
in 1971, in case of a Sturm-Liouville operator on a half-line, perturbation being
a perturbation of the boundary condition, so that in this case $V$ was one-dimensional.
An important contribution to spectral averaging was made by A.\,B.\,Alexandrov \cite{Al87}.
In 1998, B.\,Simon \cite[Theorem 1]{Si98PAMS} gave a simple short proof of the Birman-Solomyak formula.
He also noticed, that this formula holds for the wide class of Schr\"odinger operators on $\mbR^n$ \cite[Theorems 3,4]{Si98PAMS}.
The connection of this formula with the integral formula for spectral flow from non-commutative geometry is outlined in \cite{ACS}.
An interesting approach to spectral averaging via Herglotz functions can be found in \cite{GM03AA}.

In this note we present an alternative viewpoint to the spectral shift function,
and generalize the result of Simon so that it becomes applicable to a class of Dirac operators as well.
% We consider semifinite von Neumann algebra setting.

The new point of view, which the Birman-Solomyak formula suggests, is that there is a more fundamental
notion than that of the spectral shift function. We call this notion
the \emph{speed of spectral flow} or \emph{infinitesimal spectral flow}
of a self-adjoint operator $H$ under perturbation by a bounded self-adjoint operator $V.$
It was introduced in \cite{ACS} in the case of operators with compact resolvent.
It is defined by formula
\begin{gather} \label{F: inf SF (0)}
  \Phi_H(V)(\phi) = \Tr(V \phi(H)), \quad \phi \in C_c^\infty(\mbR),
\end{gather}
whenever this definition makes sense. This naturally leads to the notion of \emph{trace compatibility}
of two operators. We say that operators $H$ and $H+V$ are trace compatible, if for all $\phi \in C_\csupp^\infty(\mbR)$
the operator $V\phi(H)$ is trace class.
The spectral shift function between two trace compatible operators is then considered
as the integral of infinitesimal spectral flow. It turns out that the spectral shift function
does not depend on the path connecting the initial and final operators, a fact which follows from the
aforementioned result of B.\,Simon in the case of trace class perturbations.

The results of this note are summarized in Theorem \ref{T: 5}.
This theorem extends formulae (\ref{F: Krein's formula}) and (\ref{F: BS formula})
to the setting of trace compatible pairs $(H,H+V)$ and also strengthens
\cite[Theorems 3,4]{Si98PAMS} in the sense that it does not require $H$
to be a positive operator and maximally weakens conditions on the path $H+rV,$ $r \in [0,1].$
Our results also hold for a more general setting, when $H=H^*$ is affiliated with a semifinite
\vNa\ $\clN$ and $V=V^* \in \clN.$

Our investigation here also strengthens the link between the theory of
the Krein spectral shift function and that of spectral flow firstly discovered in \cite{ACDS}.
For exposition of the latter theory we refer to \cite{BCPRSW} and a detailed discussion of the connection
between the two theories in the situation where the resolvent of $H$ is $\tau$-compact
(here, $\tau$ is an arbitrary faithful normal semifinite trace on $\clN$)
is contained in \cite{ACS}. It should be pointed out here that the idea of viewing
the spectral shift function as the integral of infinitesimal spectral flow is akin
to \Singer's ICM-1974 proposal
to define the $\eta$ invariant (and hence spectral flow) as the integral of a one form.
Very general formulae of that type
have been produced in the framework of noncommutative geometry (see \cite{BCPRSW} and
references therein). We believe that our present approach will have applications to noncommutative
geometry, in particular, it may be useful in avoiding "summability constraints" on $H$ customarily used in that theory.

In semifinite von Neumann algebras $\clN$ Krein's formula (\ref{F: Krein's formula})
was proved for the first time in \cite{CP77ActM} in case of a bounded self-adjoint operator $H \in \clN$
and a trace class perturbation $V=V^* \in \LpN{1}$
and in \cite{ADS} for self-adjoint operators $H$ affiliated with $\clN.$

An additional reason to call $\Phi_H(V)$ the speed of spectral flow
is the following observation. Let $H$ be the operator of multiplication by $\lambda$
on $L^2(\mbR, d\rho(\lambda))$ with some measure $\rho$ and let the perturbation $V$
be an integral operator with a sufficiently regular (for example $C^1$) kernel $k(\lambda',\lambda).$
%Let $H$ be a Sturm-Liouville operator and let
%$U : L^2(\mbR,dx) \to L^2(\spectrum{H}, d\rho(\lambda))$ be a generalized Fourier transform,
%where $\spectrum{H}$ is the spectrum of $H$ and $\rho(\lambda)$ is the density of states of $H.$
%Let $V = v(x)$ be a perturbation and assume that $U v(x) U^{-1}$ is an integral operator on $L^2(\spectrum{H}, d\rho(\lambda))$
%with kernel $k(\lambda, \lambda').$
% \Tr(U V U^{-1} \phi(\lambda)) =
Then for any test function $\phi \in \Cc$
$$
  \Phi_H(V)(\phi) = \Tr(V\phi(H)) = \int_{\spectrum{H}} k(\lambda,\lambda)\phi(\lambda)\,d\rho(\lambda).
$$
Hence, the infinitesimal spectral flow of $H$ under perturbation by $V$ is the measure on the spectrum of $H$
with density $k(\lambda,\lambda)\,d\rho(\lambda).$ We note that this agrees with the classical
formula \cite[(38.6)]{LL3}
$$
  E_n^{(1)} = V_{nn}
$$
from formal perturbation theory. Here $E_n^{(0)}$ is the $n$-th eigenvalue of the unperturbed operator $H_0,$
$E_n^{(j)}, \ j=1,2,\ldots$ is the $j$-th correction term
for the $n$-th eigenvalue $E_n$ of the perturbed operator $H=H_0+V$
in the formal perturbation series $E_n = E_n^{(0)}+E_n^{(1)}+E_n^{(2)}+\ldots,$
and $V_{mn} = \big\la \psi_m^{(0)} | V | \psi_n^{(0)} \big\ra$
is the matrix element of the perturbation operator $V$ with respect to the eigenfunctions
$\psi_m^{(0)}$ and $\psi_n^{(0)}$ of the unperturbed operator $H_0$ \cite{LL3}.

Acknowledgement. We thank Alan Carey for useful comments and criticism.

%This note covers also the important case of Dirac operators, which are not considered in \cite{Si98PAMS}.

\section{Results} \nopagebreak
Let $\clN$ be a \vNa\ on a Hilbert space $\hilb$ with faithful normal semifinite trace~$\tau.$
Let $\clA = H_0 + \clA_0$ be an affine space of self-adjoint operators
affiliated with $\clN,$ where $H_0$ is a self-adjoint operator affiliated with
$\clN$ and $\clA_0$ is a vector subspace of the real Banach space of all
self-adjoint operators from $\clN.$ We say that $\clA$ is \emph{trace compatible}, if
for all $\phi \in \Cc ,$ $V \in \clA_0$ and $H \in \clA$
\begin{gather}  \label{F: V phi(H) in L1}
  V \phi(H) \in \LpN{1},
\end{gather}
where $\LpN{1}$ is the ideal of trace class operators from $\clN,$
%Clearly, $\clA$ is trace compatible iff (\ref{F: V phi(H) in L1}) holds for all bounded measurable $\phi$
%with compact support, or, iff (\ref{F: V phi(H) in L1}) holds for all $\phi = \chi_{[a,b]},$ $a<b.$
and if $\clA_0$ is endowed with a locally convex topology which coincides with
or is stronger than the uniform topology,
such that the map $(V_1, V_2) \in \clA_0^2 \mapsto V_1 \phi(H_0+V_2)$ is $\clL^1$ continuous for all
$H_0 \in \clA$ and $\phi \in \Cc.$ In particular, $\clA$ is a locally convex affine space.
The ideal property of $\LpN{1}$ and \cite[Theorem VIII.20(a)]{RS1} imply that, in the definition of
trace compatibility, the condition $\phi \in \Cc$ may be replaced by $\phi \in C_c(\mbR).$
It follows from the definition of the topology on $\clA_0$ that $H \in \clA \mapsto e^{itH}$ is norm continuous.

If $\clA = H_0 + \clA_0$ is a trace compatible affine space then we define
a (generalized) one-form (on $\clA$) of \emph{infinitesimal spectral flow}
or \emph{speed of spectral flow} by the formula
\begin{gather} \label{F: inf SF}
  \Phi_H(V) = \tau(V \delta(H)), \quad H \in \clA, \ V \in \clA_0,
\end{gather}
where $\delta$ is Dirac's delta function. The last formula is to be understood
in a generalized function sense, i.e.
$
  \Phi_H(V)(\phi) = \tau(V \phi(H)), \ \ \phi \in \Cc.
$
$\Phi$ is a generalized function, since if $\phi_n \to 0$ in $\Cc $
such that $\supp(\phi_n) \subseteq \Delta$ then
$\abs{\tau(V \phi_n(H))} \leq \norm{V E_\Delta^H}_1 \norm{\phi_n(H)} \to 0.$
Here $\norm{A}_1 = \tau(|A|).$

Since $\phi$ can be taken from $C_c(\mbR),$ for each $V \in \clA_0$ the infinitesimal spectral flow
$\Phi_H(V)$ is actually a measure on the spectrum of $H.$

By a smooth path $\set{H_r}_{r \in \mbR}$ in $\clA,$ we mean a differentiable path, such that its
derivative $\frac {dH_r}{dr} \in \clA_0$ is continuous.

Let $\Pi = \set{(s_0,s_1) \in \mbR^2 \colon s_0s_1 \geq 0, \abs{s_1} \leq \abs{s_0}},$
and let $$d\nu_f(s_0,s_1) = \sgn(s_0)\frac i{\sqrt{2\pi}} \hat f(s_0)\,ds_0\,ds_1.$$
If $f \in C^2_\csupp(\mbR)$ then $(\Pi,\nu_f)$ is a finite measure space \cite{ACDS}.
For any $H_0,H_1 \in \clA,$ any $X \in \clA_0$ and any non-negative $f \in \Cc$ set by definition
\begin{multline} \label{F: def of TH1H0(X)}
    \TT^{H_1,H_0}_{f^{[1]}}(X) = \int_\Pi \big( e^{i(s_0-s_1)H_1} \sqrt{f}(H_1) X e^{i s_1 H_0}
    \\ + e^{i(s_0-s_1)H_1} X \sqrt{f}(H_0) e^{i s_1 H_0}\big)\,d\nu_{\sqrt{f}} (s_0,s_1),
\end{multline}
where the integral is taken in the $so^*$-topology.
For justification of this notation and details see \cite{ACS}.
\begin{lemma} \label{L: 1} If $\set{H_r} \subset \clA$ is a path,
continuous (smooth) in the topology of $\clA_0,$ and if $f \in C_c^2(\mbR)$ then
$$
  r \mapsto f(H_r) - f(H_0)
$$
takes values in $\LpN{1}$ and it is $\LpN{1}$ continuous (smooth).
\end{lemma}
\begin{proof}
We can assume that $f$ is non-negative and that $\sqrt f \in C_c^2(\mbR).$
It is proved in \cite{ACS} that
\begin{gather} \label{F: phi Hr - phi H0}
  f(H_r) - f(H_0) = \TT^{H_r,H_0}_{f^{[1]}}(H_r-H_0).
\end{gather}
Since $e^{i(s_0-s_1)x} \sqrt f(x),$ $e^{is_1x} \sqrt f(x) \in C_c^2(\mbR),$
trace compatibility implies that the integrand of the right hand side of
(\ref{F: def of TH1H0(X)}) takes values in $\clL^1$ and is $\clL^1$-continuous
(smooth), so the dominated convergence theorem completes the proof.
\end{proof}

If $\Gamma = \set{H_r}_{r \in [0,1]}$ is a smooth path in $\clA,$ then we define the spectral shift function $\xi$
along this path as the integral of infinitesimal spectral flow:
$
  \xi = \int_\Gamma \Phi,
$
or
\begin{gather} \label{F: def of ssf}
  \xi(\phi) = \int_0^1 \taubrs{ \frac {dH_r}{dr}\,  \phi(H_r)}\,dr, \quad \phi \in \Cc.
\end{gather}
Now we prove that the spectral shift function is well-defined
in the sense that it does not depend on the path of integration.

A one-form $\alpha_H(V)$ on an affine space $\clA$ is called \emph{exact} if there exists a zero-form $\theta_H$ on $\clA$
such that $d\theta = \alpha,$ i.e.
$$
  \alpha_H(V) = \frac d{dr} \theta_{H+rV} \Big|_{r=0}.
$$
We say that the generalized one-form $\Phi$ is exact if $\Phi(\phi)$ is an exact form for any $\phi \in \Cc.$

The proof of the following proposition follows the lines of the proof of \cite[Proposition 3.5]{ACS}.
\begin{prop}\label{P: for H0+Lp1 Phi exact}
  The infinitesimal spectral flow $\Phi$ is exact.
\end{prop}
\begin{proof}
Let $V \in \clA_0,$ $H_r = H_0+rV,$ $r \in [0,1],$ and let $f \in \Cc.$
By (\ref{F: phi Hr - phi H0})
\begin{multline}
  f(H_r) - f(H_0) = \TT^{H_r,H_0}_{f^{[1]}}(r V)
    \\ = \int_\Pi \big( e^{i(s_0-s_1)H_r} \sqrt{f}(H_r) r V e^{i s_1 H_0}
     + e^{i(s_0-s_1)H_r} r V \sqrt{f}(H_0) e^{i s_1 H_0}\big)\,d\nu_{\sqrt{f}} (s_0,s_1),
    \\ = \int_\Pi \big( e^{i(s_0-s_1)H_0} \sqrt{f}(H_0) r V e^{i s_1 H_0}
     + e^{i(s_0-s_1)H_0} r V \sqrt{f}(H_0) e^{i s_1 H_0}\big)\,d\nu_{\sqrt{f}} (s_0,s_1)
     \\ + \int_\Pi \big( (e^{i(s_0-s_1)H_r} \sqrt{f}(H_r) - e^{i(s_0-s_1)H_0} \sqrt{f}(H_0)) r V e^{i s_1 H_0}
     \\ + (e^{i(s_0-s_1)H_r} - e^{i(s_0-s_1)H_0}) r V \sqrt{f}(H_0) e^{i s_1 H_0}\big)\,d\nu_{\sqrt{f}} (s_0,s_1)
    \\ =: \TT^{\HH_0,\HH_0}_{f^{[1]}}(r\VV) +  R_1 + R_2.
\end{multline}
All three summands here are trace class by the trace compatibility assumption.
So, for any $S \in \clN$
$$
  \tau(S(f(H_r) - f(H_0))) = r \tau\brs{S \TT^{\HH_0,\HH_0}_{f^{[1]}}(\VV)} + \tau(SR_1) + \tau(SR_2).
$$
Now, Duhamel's formula and (\ref{F: phi Hr - phi H0}) show that $\tau(SR_1) = o(r)$ and $\tau(SR_2) = o(r).$
Hence,
\begin{gather} \label{F: dot h = tau tHH}
  \frac d{dr} \tau(S(f(H_r) - f(H_0))) = \tau\brs{S \TT^{H_r,H_r}_{f^{[1]}}(V)}.
\end{gather}
This implies that for any $S \in \clN$
$$
  \tau(S(f(H_1) - f(H_0))) = \tau\brs{\int_0^1 S \TT^{H_r,H_r}_{f^{[1]}}(V)\,dr}.
$$
Now let $H_0 \in \clA$ be a fixed operator and for any $f \in \Cc$ let
$$
  \theta_H^f := \int_0^1 \tau(Vf(H_r))\,dr,
$$
where $H_r = H_0 + rV,$ $H= H_1.$
We are going to show that
$d\theta_H^f(X) = \Phi_H(X)(f)$ for any $X \in \clA_0.$

Following the proof of \cite[Proposition 3.5]{ACS} we have
\begin{multline*}
  (A) := d\theta_H^f(X)
        = \lims{s\to 0} \int_0^1 \tau\Brs{X f(\HH_r + s r X)}\,dr
    \\    + \lims{s\to 0} \frac 1s \int_0^1 \tau\BRS{\VV \Brs{f(\HH_r+s r X) - f(\HH_r)} }\,dr.
\end{multline*}
By definition of $\clA_0$ topology the integrand of the first summand is continuous with respect to $r$ and $s.$
So, the first summand is equal to
$$
  \int_0^1 \tau\Brs{X f(\HH_r)}\,dr.
$$
By \cite[Theorem 5.3]{ACDS} the second summand is equal to
  \begin{multline*}
     \lims{s\to 0} \frac 1s \int_0^1 \taubrs{\VV \TT^{\HH_r+srX,\HH_r}_{f^{[1]}}(s r X)}\,dr
       = \lims{s\to 0} \int_0^1 \taubrs{\VV \TT^{\HH_r+srX,\HH_r}_{f^{[1]}}(r X)}\,dr
      \\ = \int_0^1 \taubrs{\VV\, \TT^{\HH_r,\HH_r}_{f^{[1]}}(r X)}\,dr
       = \int_0^1 \taubrs{X\, \TT^{\HH_r,\HH_r}_{f^{[1]}}(\VV)}r\,dr,
  \end{multline*}
where the second equality follows from the definition of $\clA_0$-topology
and the last equality follows from \cite[Lemma 3.2]{ACS}.
Now, using (\ref{F: dot h = tau tHH}) and integrating by parts we get
\begin{multline}
  (A) = \tau\brs{X f(H_1)-X f(H_0)} + \int_0^1 \brs{\tau(X f(\HH_r)) - \tau(X [f(\HH_r)-f(\HH_0)])}\,dr
   \\ = \tau\brs{X f(H_1)}.
\end{multline}
\end{proof}
The argument before \cite[Proposition 1.5]{CP98CJM} now implies
\begin{cor}
  The spectral shift function given by (\ref{F: def of ssf}) is well-defined.
\end{cor}

%
% Proposition \ref{P: for H0+Lp1 Phi exact} shows that in the case when $\clA_0 = \LpN{1}$
% equipped with the norm $\norm{\cdot}_{\LpN{1}}$
% the spectral shift function is well-defined.
%
\begin{prop} \label{P: 2}
If $r \in \mbR \mapsto H_r \in \clA$ is smooth then the equality
\begin{gather}\label{F: tr ddr f(Hr) = ...}
  \taubrs{\frac {d f(H_r)}{dr} \phi(H_r)} = \taubrs{\dot H_r f'(H_r) \phi(H_r)}
\end{gather}
holds for any $f \in C_c^2(\mbR)$ and any bounded measurable function $\phi.$
\end{prop}
\begin{proof} Without loss of generality, we can assume that $f \geq 0$ and $\sqrt{f} \in C_c^2(\mbR).$
We prove the above equality at $r = 0.$
%We have
%\begin{gather*}
%  \taubrs{\frac {d f(H_r)}{dr} \Big|_{r=0} \phi(H_0)}
%    = \taubrs{\phi(H_0) \frac {d (f(H_r) - f(H_0))}{dr}\Big|_{r=0}}.
%\end{gather*}
The formula (\ref{F: phi Hr - phi H0}) and the dominated convergence theorem imply that
\begin{multline*}
  \taubrs{\frac {d f(H_r)}{dr}\Big|_{r=0} \phi(H_0)}
    = \tau \Big( \phi(H_0)
    \int_\Pi \lim_{r \to 0} \big[ e^{i(s_0-s_1)H_r} \sqrt{f}(H_r) \frac{H_r-H_0}{r} e^{i s_1 H_0}
    \\ + e^{i(s_0-s_1)H_r} \frac{H_r-H_0}{r} \sqrt{f}(H_0) e^{i s_1 H_0}\big]\,d\nu_{\sqrt{f}} (s_0,s_1)
   \Big),
\end{multline*}
where the limit is taken in $\LpN{1}.$ By the $\clA_0$-smoothness of $\set{H_r},$ we have
\begin{multline*}
  \taubrs{\frac {d f(H_r)}{dr}\Big|_{r=0} \phi(H_0)}
    = \tau \Big( \phi(H_0)
    \int_\Pi \big[ e^{i(s_0-s_1)H_0} \sqrt{f}(H_0) \dot H_r \big|_{r=0} e^{i s_1 H_0}
    \\ + e^{i(s_0-s_1)H_0} \dot H_r\big|_{r=0} \sqrt{f}(H_0) e^{i s_1 H_0}\big]\,d\nu_{\sqrt{f}} (s_0,s_1)
   \Big),
\end{multline*}
so that by \cite[Lemmas 3.7, 3.10]{ACDS} and letting $A = \dot H_r \big|_{r=0}$
\begin{multline*}
  \taubrs{\frac {d f(H_r)}{dr}\Big|_{r=0} \phi(H_0)}
    = 2 \int_\Pi \tau \Big( \phi(H_0)
    e^{i s_0 H_0} \sqrt{f}(H_0) A \Big)\,d\nu_{\sqrt{f}} (s_0,s_1)
    \\ = 2 \tau \Big( A  \phi(H_0) \sqrt{f}(H_0) \int_{-\infty}^\infty
    e^{i s_0 H_0} \frac i{\sqrt{2\pi}} s_0 \fourier {\sqrt{f}}(s_0)\,ds_0\Big)
     = \tau \big( A \phi(H_0) f'(H_0) \big).
\end{multline*}
%It can be proved by using (\ref{F: phi Hr - phi H0}).
\end{proof}

\begin{prop} \label{T: 4}
The spectral shift function given by (\ref{F: def of ssf}) satisfies Krein's formula, i.e. for any $f \in C^\infty_c,$
$H_0, H_1 \in \clA$
$$
  \taubrs{f(H_1) - f(H_0)} = \xi(f').
$$
\end{prop}
\begin{proof}
Taking the integral of (\ref{F: tr ddr f(Hr) = ...}) with $\phi = 1$
we have
$$
  \int_0^1 \taubrs{\frac {d (f(H_r) - f(H_0))}{dr}}\,dr = \int_0^1 \taubrs{\dot H_r f'(H_r)}\,dr.
$$
The right hand side is $\xi(f')$ by definition.
It follows from Lemma \ref{L: 1} that one can interchange the trace and the derivative
in the left hand side.
\end{proof}
\begin{cor} In case of trace class perturbations,
the spectral shift function $\xi$ defined by (\ref{F: def of ssf})
coincides with classical definition, given by \cite[Theorem 3.1]{ADS}.
\end{cor}
\begin{proof} This follows from Theorem 6.3 and Corollary 6.4 of \cite{ACDS}.
\end{proof}

For trace class perturbations
the absolute continuity of the spectral shift function is established in \cite{Kr53MS} (see also \cite{GM03AA}).
For the general semifinite case we refer to \cite{ADS,ACDS}.

\begin{lemma} \label{L: some equality}
Let $\clA$ be a trace compatible affine space and let $f \in \Cc.$
Let $H_0, H_1 \in \clA,$ let $\xi$ and $\xi_f$
be the spectral shift distributions of the pairs $(H_0,H_1)$
and $f(H_0),f(H_1)$ respectively.
Then for any $\phi \in \Cc$
  \begin{gather} \label{F: xi f = xi ...}
    \xi_f(\phi) = \xi(\phi\circ f \cdot f').
  \end{gather}
\end{lemma}
\begin{proof} By Proposition \ref{P: 2} for any $f, \phi \in \Cc$
\begin{gather*}
  \taubrs{ \frac {d f(H_r)}{dr}  \phi(f(H_r))} = \taubrs{ \dot H_r f'(H_r) \phi(f(H_r))}
  = \taubrs{ \dot H_r (F \circ f)'(H_r)},
\end{gather*}
where $F' = \phi.$ Hence, for any smooth path $\Gamma = \set{H_r}_{r \in [0,1]} \subseteq \clA$
\begin{gather} \label{F: int G = int G}
  \int_0^1 \taubrs{ \frac {d f(H_r)}{dr}  \phi(f(H_r))}\,dr = \int_0^1 \taubrs{ \dot H_r (F \circ f)'(H_r)}\,dr,
\end{gather}
which is (\ref{F: xi f = xi ...}).
% Since the path $f(H_r)$ lies in a standard trace compatible affine space $H_0+\LpN{1},$
% the left hand side integral does not depend on a path by Proposition \ref{P: for H0+Lp1 Phi exact},
% hence so does the right hand side.
% Finally, any $\psi \in \Cc$ can be written as $\psi = (F \circ f)'$
% with $F', f \in \Cc.$  %%% <--- it is a mistake
\end{proof}

\begin{prop} Let $\clA = H_0 +\clA_0$ be a trace compatible affine space and
let $\clA_0$ be such that for any $V \in \clA_0$ there exist positive $V_1, V_2 \in \clA$
such that $V = V_1 - V_2.$ Then the spectral shift function $\xi$ of any pair $H,H+V \in \clA$
is absolutely continuous.
\end{prop}
\begin{proof}
Since the map $(V_1,V_2) \mapsto \taubrs{V_1 \phi(H+V_2)}$ is $\LpN{1}$\tire continuous
(by definition), it follows that the infinitesimal spectral flow is a uniformly locally finite measure
with respect to the path parameter.
Hence, the spectral shift function is also a locally finite measure being the integral
of locally finite measures, which are uniformly bounded on every segment.

If, for $H_0, H_1, H_2 \in \clA,$ the spectral shift functions from $H_0$ to $H_1$ and from $H_1$ to $H_2$
are absolutely continuous, then evidently the spectral shift function from $H_0$ to $H_2$ is also absolutely
continuous.
Hence, if $V = V_1 - V_2$ with $0 \leq V_1 , V_2 \in \clA_0,$
then representing the spectral shift function from $H$ to $H+V$ as the sum of spectral shift function
from $H$ to $H+V_1$ and from $H+V_1$ to $H+V,$ we see that we can assume that the perturbation $V$
is positive.

By Lemma \ref{L: 1} and \cite[Theorem 3.1]{ADS}
the spectral shift function $\xi_f$ of the pair $(f(H), f(H+V))$
is absolutely continuous.
Let us suppose that the spectral shift function $\xi$ of the pair $(H,H+V)$
has non-absolutely continuous part $\mu.$

Without loss of generality, we can assume that
there exists a set of Lebesgue measure zero $E \subset (\eps, 1-\eps)$ such that $\mu(E) > 0.$
For any $a, b \in \mbR$  with $b-a >2$ let us consider a "cap"-function $f_{a,b},$ i.e. $f$ is a smooth function which is zero
on $(-\infty,a)$ and $(b,\infty),$ it is $1$ on $(a+1, b-1)$
and its derivatives on $(a+\eps,a+1-\eps)$ and $(b-1+\eps,b-\eps)$ is $1$ and $-1$ respectively.

Let $U$ be an open set of Lebesgue measure $ < \delta$ such that $E \subset U$
and let $\phi$ be a smoothed indicator of $U.$
Then (\ref{F: xi f = xi ...}),
applied to functions $\phi$ and $f_{0,b}$ and to functions
$\phi$ and $f_{a,b}$ (with big enough $b$) implies that
$\mu(E) = \mu(a+E),$ i.e. $\mu$ is translation invariant. Since it is also locally finite
it is some multiple of Lebesgue measure. This yields a contradiction.
\end{proof}

We summarize the results in the following theorem.
\begin{thm} \label{T: 5} Let $\clA$ be a trace compatible affine space of operators in
a semifinite von Neumann algebra $\clN$ with a normal semifinite faithful trace~$\tau.$
Let $H$ and $H+V$ be two operators from $\clA.$ Let the spectral shift
(generalized) function $\xi_{H,H+V}$
be defined as the integral of infinitesimal spectral flow by the formula
$$
  \xi_{H,H+V}(\phi) = \int_\Gamma \Phi(\phi) = \int_0^1 \Phi_{H_r}(\dot H_r)(\phi)\,dr, \quad \phi \in C^\infty_\csupp,
$$
where $\Gamma = \set{H_r}_{r\in [0,1]}$ is any piecewise smooth path in $\clA$ connecting $H$ and $H+V.$
Then the spectral shift function is well-defined in the sense
that the integral does not depend on the choice of the
piecewise smooth path $\Gamma$ connecting $H$ and $H+V,$
and it satisfies Krein's formula
$$
  \tau(f(H+V) - f(H)) = \xi(f'), \quad f \in C^\infty_\csupp.
$$
Moreover, if for any $V \in \clA_0$ there exist $V_1, V_2 \in \clA_0$ such that $V = V_1 - V_2,$
then $\xi_{H,H+V}$ is an absolutely continuous measure.
\end{thm}

%\vskip 0.2 cm
\medskip

Two extreme examples of trace compatible affine spaces are
$H_0+\saLpN{1},$ $H_0=H_0^* \eta \clN,$ with the topology induced by $\LpN{1}$ \cite{ACDS,ADS},
and $D_0+\clN_{sa},$ where $(D_0-i)^{-1}$ is $\tau$\tire compact,
with the topology induced by operator norm \cite{ACS}. In particular,
the space $-\Delta + C(M),$ where $(M,g)$ is a compact Riemannian manifold,
and $\Delta$ is the Laplacian, is trace compatible.

As an example of an intermediate trace compatible affine space one can consider Schr\"odinger operators $-\Delta + C_c(\mbR^n)$
with the inductive topology of uniform convergence.
It is proved in \cite[Section B9]{Si82BAMS} that for this example
the condition (\ref{F: V phi(H) in L1}) holds. It also follows from \cite[Section B9]{Si82BAMS}
that $\norm{g f(H)}_1 \leq C \norm{g}_2,$ where $C$ depends only on $f,$ on the support of $g$ and on $\norm{V_-}_\infty,$
where $V_-$ is the negative part of $V,$ $H = -\Delta + V.$ So, the condition on the topology of $\clA_0$ is fulfilled
by~(\ref{F: phi Hr - phi H0}).

Another example is given by Dirac operators of the form $D + \clA_0,$
where $D = \sum\limits_{j=1}^n \alpha_j \frac \partial{\partial x_j},$
$\alpha_1,\ldots,\alpha_n$ are $m\times m$-matrices such that
$\alpha_j\alpha_k + \alpha_k\alpha_j = -2\delta_{jk},$ and
$$
  \clA_0 = \set{ a = a^* \in C_c(\mbR^n, \mathbf{M}_m(\mbR)) \colon \exists \phi = \phi^* \in C^1(\mbR^n) \ \ iD\phi = a}
$$
with the inductive topology of uniform convergence.
A proof that the space $D + \clA_0$ is trace compatible
can be reduced to \cite[Theorem 4.5]{SimTrId} via the gauge transformation $\psi \mapsto e^{-i\phi(x)} \psi.$
We have
$$
  (D+a)(e^{-i\phi(x)} u) = \sum\limits_{j=1}^n\brs{-i e^{-i\phi(x)}\frac \partial{\partial x_j} \phi(x) \alpha_j u
  + e^{-i\phi(x)} \alpha_j  \frac \partial{\partial x_j} u} + a e^{-i\phi(x)} u,
$$
where $u$ is an $m$-column of $C^\infty_\csupp$-functions.
So, if $iD\phi = a$ then
$$
  e^{i\phi(x)}(D+a)(e^{-i\phi(x)} u) = Du.
$$
Hence, $(D + a)^2 = e^{-i\phi(x)} D^2 e^{i\phi(x)}.$
This shows that $g f((D + a)^2),$ $g, a \in \clA_0, \ f \in \Cc,$ is trace class iff $gf(e^{-i\phi(x)}D^2 e^{i\phi(x)})
=
e^{-i\phi(x)}g e^{i\phi(x)}f(e^{-i\phi(x)}D^2 e^{i\phi(x)})$ is trace class.
%where $\phi$ is a solution of the equation $a = iD\phi.$
But the last operator is unitarily equivalent to
$g f(D^2),$ which is trace class by \cite[Theorem 4.5]{SimTrId}. So,
if $f \geq 0$ then by the same theorem
\begin{gather} \label{F: 1-norm estimate from Simon}
  \norm{g f(D + a)}_1 \leq C \norm{g}_\infty \norm{f}_\infty,
\end{gather}
where $C$ depends on supports of $g$ and $f.$ Hence, for $g,g_1,a,a_1 \in \clA_0,$
we have
\begin{gather*}
  \norm{gf(D+a) - g_1f(D+a_1)}_1 \leq \norm{(g-g_1)f(D+a)}_1 + \norm{g_1(f(D+a) - f(D+a_1)}_1.
  %\\ \leq C \norm{g-g_1}_\infty \norm{f}_\infty + ...,
\end{gather*}
So, the condition on the topology of $\clA_0$ is fulfilled by (\ref{F: phi Hr - phi H0}) and
(\ref{F: 1-norm estimate from Simon}).

In case $n=m=1,$ we have $D + \clA_0 = \frac 1i \frac d{dx} + C_c(\mbR).$

\mathsurround 0pt
\ndef{\AndSoOn}{$\dots$}

%\end{document}
%  End of MyListOfRef.tex
%\input ../ListOfRef/MyListOfRef
\end{document}